\RequirePackage{fix-cm}

\documentclass{svjour3}

\newif\ifTikz
\Tikztrue

\smartqed

\usepackage{graphicx}
\usepackage{amsmath,amsfonts}
\usepackage{tensor}
\usepackage{subfig}
\usepackage[rightcaption]{sidecap}

\usepackage[colorlinks=true,citecolor=blue,urlcolor=blue,linkcolor=black]{hyperref}

\newcommand{\xx}{\mathbf{x}}

\begin{document}

\title{A New Family of Regularized Kernels for the Harmonic Oscillator}

\author{
  Benjamin W.~Ong \and
  Andrew J.~Christlieb \and
  Bryan D.~Quaife
}
\authorrunning{B.~W.~Ong, A.~J.~Christlieb, B.~D.~Quaife}

\institute{
  Benjamin W.~Ong \at
  Mathematical Sciences, Michigan Technological University, Houghton, MI, 49931\\
  {\tt ongbw@mtu.edu}
  \and
  Andrew J.~Christlieb \at Dept.~of
  Mathematics, Michigan State University, East Lansing, MI, 48824 \\
  {\tt andrewch@msu.edu}
  \and
  Bryan D.~Quaife \at Department of Scientific Computing, Florida
    State University, Tallahassee, FL, 32306 \\{\tt bquaife@fsu.edu}
}

\date{Received: \today}

\maketitle

\begin{abstract}
  In this paper, a new two-parameter family of regularized kernels is
  introduced, suitable for applying high-order time stepping to N-body
  systems.  These high-order kernels are derived by truncating a
  Taylor expansion of the non-regularized kernel about
  $(r^2+\epsilon^2)$, generating a sequence of increasingly more
  accurate kernels.  This paper proves the validity of this
  two-parameter family of regularized kernels, constructs error
  estimates, and illustrates the benefits of using high-order kernels
  through numerical experiments.
\end{abstract}

\section{Introduction}

\subsection{Problem statement}
This paper is concerned with solutions to the system
\begin{align}
  \ddot{\xx}_{j} = -\sum_{k \neq j} w_{jk} \nabla G(\xx_{j} - \xx_{k}),
  \label{eqn:hamiltonian_system}
\end{align}
where $w_{jk}$ are constants, and $G$ is the fundamental solution to
Laplace's equation. 
Specifically,
\begin{align*}
  G(r) = 
  \begin{cases}
    \vspace*{0.2cm}
    \displaystyle
    -\frac{r}{2} & \text{in } \mathbb{R}^1,\\
    \vspace*{0.2cm}
    \displaystyle
    -\frac{\ln{r}}{2\pi} & \text{in } \mathbb{R}^2,\\
    \displaystyle
    \frac{1}{4\pi r} & \text{in } \mathbb{R}^3,
  \end{cases}
\end{align*}
where $r = \|\xx\|_{2}$.  Henceforth, $\| \cdot \|_{2}$ is denoted as
$\| \cdot \|$.  System~\eqref{eqn:hamiltonian_system}, modulo a change
in sign, arises in many dynamical systems such as the dynamics of
charged particles~\cite{jackson1999}, vortex
dynamics~\cite{leo1980,win-leo1993}, and planetary
motions~\cite{bate1971}.  The weights $w_{jk}$ can be interpreted as
the interaction between point masses at $x_j$ and $x_k$, or viewed as
quadrature weights that approximate interactions due to a distribution
of masses.  By taking the first integral of
system~\eqref{eqn:hamiltonian_system}, the Hamiltonian is obtained,
\begin{align*}
  H(\xx,\dot{\xx}) = \frac{1}{2}\sum_{j} \|\dot{\xx}_{j}\|^{2} + 
    \sum_{j}\sum_{k \neq j} w_{jk} G(\xx_{j} - \xx_{k}),
\end{align*}
which is a conserved quantity.  

\begin{figure}[htbp]
  \begin{minipage}{0.48\textwidth}
    \includegraphics[scale=1]{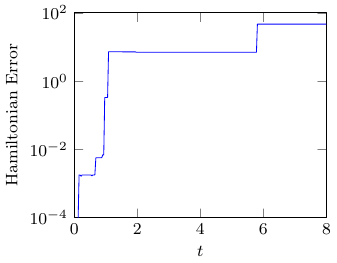}
  \end{minipage}
  \quad
  \begin{minipage}{0.48\textwidth}
    \includegraphics[scale=1]{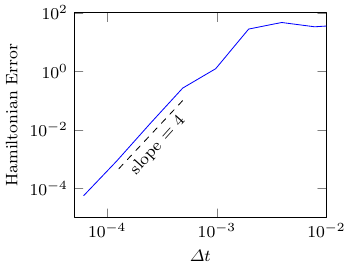}
  \end{minipage}
  \caption{\label{fig:Hamiltonian_eps0} {\em Left}: The Hamiltonian as
    a function of time when equation~\eqref{eqn:hamiltonian_system} is
    integrated with a fourth-order symplectic integrator with $\Delta
    t = 3.91 \times 10^{-3}$.  {\em Right}: The error of the
    Hamiltonian at the time horizon $T=8$ for various time step sizes.
    Fourth-order convergence is eventually attained for small time
    steps.  The dashed black line indicates fourth-order convergence.}
\end{figure}

In numerical simulations, controlling the error in the Hamiltonian is
of the utmost importance.  Special classes of integrators, such as
symplectic integrators~\cite{forest1990fourth,4332919} or
energy-conserving integrators~\cite{faou2004}, have been designed to
help preserve either the symplectic structure of the equations
(thereby adding stability and controlling the Hamiltonian error in
some fashion as the simulation progresses), or explicitly conserving
the Hamiltonian.  If a fourth-order symplectic
integrator~\cite{forest1990fourth} with a modest time step is used to
solve a 25-body problem in $\mathbb{R}^{2}$ with $w_{jk} = \pm 1$, the
error in the Hamiltonian rises quickly; in fact, each rise or dip
corresponds to when particles cross each other.  This is illustrated
in the left plot of Figure~\ref{fig:Hamiltonian_eps0}.  A convergence
study shows poor convergence to the true Hamiltonian for larger
$\Delta t$ before fourth-order convergence is observed for
sufficiently small $\Delta t$.  The results of the convergence study
is illustrated in the right plot of Figure~\ref{fig:Hamiltonian_eps0}.

\begin{figure}[htbp]
  \begin{minipage}{0.48\textwidth}
    \includegraphics[scale=1]{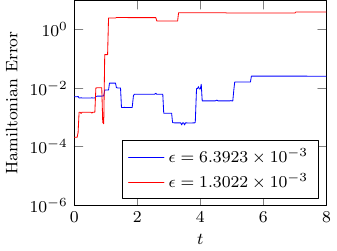}
  \end{minipage}
  \quad
  \begin{minipage}{0.48\textwidth}
    \includegraphics[scale=1]{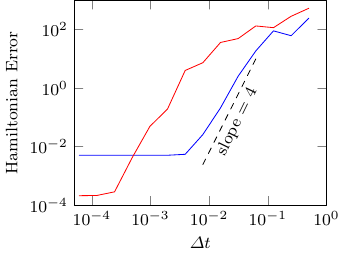}
  \end{minipage}
  \caption{\label{fig:Hamiltonian_eps1} {\em Left}: The Hamiltonian as
  a function of time when the regularized
  equation~\eqref{eqn:regularized_system} is integrated
  with a fourth-order symplectic integrator with $\Delta t = 3.91
  \times 10^{-3}$.  Although a smaller $\epsilon$ results in a smaller
  modelling error initially (at $t=0$), the steeper gradient in the
  kernel causes the time discretization error to dominate.  {\em
  Right}: The error of the Hamiltonian at the time horizon $T=8$ for
  various time step sizes.  The dashed black line indicates
  fourth-order convergence.}
\end{figure}

A common way to overcome the reduced order of accuracy is to solve a
regularized system instead of the original system,
\begin{align}
  \ddot{\xx}_{j} = -\sum_{k \neq j} w_{jk}
    \nabla G^{\epsilon}(\xx_{j} - \xx_{k}),
  \label{eqn:regularized_system}
\end{align}
which has the modified Hamiltonian
\begin{align*}
  H^{\epsilon}(\xx,\dot{\xx}) = \frac{1}{2}\sum_{j} \|\dot{\xx}_{j}\|^{2} + 
    \sum_{j}\sum_{k \neq j} w_{jk} G^{\epsilon}(\xx_{j} - \xx_{k}).
\end{align*}
The following algebraic regularization have been used in
plasma and vortex simulations~\cite{christlieb:gfp06,lindsay-krasny01},
\begin{align}
  G^{\epsilon}(r) = 
  \begin{cases}
    \displaystyle
    -\frac{1}{2}\sqrt{r^2+\epsilon^2}
    & \text{in }\mathbb{R}^1, \\
    \displaystyle
    -\frac{\ln{\sqrt{r^2+\epsilon^2}}}{2\pi}
    & \text{in }\mathbb{R}^2, \\
    \displaystyle
    \frac{1}{4\pi}\frac{1}{\sqrt{r^2+\epsilon^2}}
    & \text{in }\mathbb{R}^3.
  \end{cases}
  \label{eqn:regularized_n0}
\end{align}
We shall refer to this algebraic regularization as the one-parameter
family of regularized kernels.  In Figure~\ref{fig:Hamiltonian_eps1},
the same fourth-order symplectic integrator is used to solve
system~\eqref{eqn:regularized_system} for the same 25-body system
described earlier. Two different values of $\epsilon$ are used to
specify the regularized kernel.  The regularized system achieves
fourth-order accuracy for large time step sizes.  However for small
time steps, the error stagnates; the value at which the error
stagnates corresponds to the difference between the unregularized and
regularized Hamiltonians, $|H(0)-H^{\epsilon}(0)|$.  We shall refer to
this difference as modelling error.  This modelling error is observed
in the left plot of Figure~\ref{fig:Hamiltonian_eps1} at $t=0$.

The modelling error can be reduced by decreasing $\epsilon$.  However,
this results in steeper derivatives of $G^{\epsilon}(r)$, creating
larger jumps in the Hamiltonian error (for large time steps) as two
particles pass one another.  The left plot of
Figure~\ref{fig:Hamiltonian_eps1} demonstrates this behavior.  The net
effect, as can be seen in the right plot of
Figure~\ref{fig:Hamiltonian_eps1}, is that smaller values of
$\epsilon$ require smaller time step sizes before the smaller modelling error is realized.
This paper constructs a two-parameter family of regularized kernels
that allows an integrator to achieve small modelling error with large
time steps.

\subsection{Related work}
One approach for forming a regularization is to solve $\mathcal{L}
G^{\epsilon} = \delta_{\epsilon}$, where $\delta_{\epsilon}$ is a
regularized approximation of the delta distribution, sometimes
referred to as blob approximations when studying vortex
dynamics~\cite{cortez2001,Cortez2000428}.  This approach has also been
applied in other fields, for example in plasma
physics~\cite{majda-PD94}, where a piecewise approximation to the
delta function was introduced.  Although the blob approximations
converge to the unregularized system as $\epsilon\to0$, blobs with
small $\epsilon$ generate large variations in $\frac{\partial
  G^{\epsilon}}{\partial r}$, leading to severe time-step restrictions
when resolving the evolution dynamics.  The quality of different delta
distribution regularizations has recently been analyzed in a
functional analysis setting~\cite{hos-nig-sto2016}.  Alternatively,
one can replace the fundamental solution with a regularized version.
For example, a regularized fundamental solution might satisfy $\Delta
G^\epsilon(r) = n(\epsilon)\Delta G(r)
s(r/\epsilon)$~\cite{beale2004gbb}, where $s(r)$ is a shape function
satisfying
\begin{align*}
  s(r) = \text{erf}(r) - \frac{2}{\sqrt{\pi}}re^{-r^2},
\end{align*}
sometimes referred to as a Gaussian mollifer, and $n(\epsilon)$ is
some normalizing factor that depends on the regularization $\epsilon$.
Different filters or mollifiers can also be used~\cite{MR3101517}, but
as before, small $\epsilon$ generate large variations in
$\frac{\partial G^{\epsilon}}{\partial r}$.  This paper replaces the
fundamental solution with a two-parameter algebraic regularized
kernel. Previously, one-parameter algebraic regularizations,
equation~\eqref{eqn:regularized_n0}, have been used in in plasma
physics~\cite{christlieb:gfp06,gibbon2010,win-leo1993}, and for vortex
sheet computations in fluid
dynamics~\cite{leo1980,lindsay-krasny01,MR2520287}.  The modeling
error that arises from the one-parameter family of regularized kernels
is undesirable for practical values of $\epsilon$.

\subsection{Paper outline}
This paper seeks new regularizations that reduce modelling error while
allowing an integrator to realize these improved modelling error with
large time steps.  In Section~\ref{sec:regularized_potentials}, a
two-parameter family of regularized kernels is introduced along with a
notion of a global smoothing error, which will be used to quantify the
quality of the regularized kernels and how they impact a numerical
simulation.  In Section~\ref{sec:validity}, the Laplacian of the
two-parameter family of regularized potentials is shown to converge to
the delta function.  In Section~\ref{sec:analysis}, the global
smoothing errors of the two-parameter family of regularized kernels is
analyzed.  The regularized potentials are used to solve various
$N$-body problems in Section~\ref{sec:numerics}.  Finally, in
Section~\ref{sec:conclusions}, we summarize the results and discuss
future work.

\section{Regularized kernels}
\label{sec:regularized_potentials}

\subsection{Two-parameter family of regularized kernels}
A two-parameter family of algebraic regularized kernels can be
constructed by taking the binomial or Taylor expansion of the
non-regularized kernels,
\begin{align*}
G(r) =
\begin{cases}
-\displaystyle\frac{1}{2}(r^2+\epsilon^2 -\epsilon^2)^{1/2}
& \text{in }\mathbb{R}^1, \\
\displaystyle
-\frac{1}{4\pi} \ln(r^{2} + \epsilon^{2}) - 
\displaystyle
\frac{1}{4\pi}\ln(1-\epsilon^{2}(r^{2}+\epsilon^{2})^{-1})
& \text{in }\mathbb{R}^2, \\
\displaystyle
\frac{1}{4\pi}(r^2+\epsilon^2 -\epsilon^2)^{-1/2}
& \text{in }\mathbb{R}^3,
\end{cases}
\end{align*}
and then truncating the expansion after $n$ terms,
\begin{align*}
  G^{\epsilon,n}(r) = 
\begin{cases}
-\displaystyle\frac{1}{2}\sum_{\ell=0}^n {\frac{1}{2}
    \choose \ell}\left(-\epsilon^2 \right )^\ell \left (r^2+\epsilon^2\right
  )^{1/2-\ell}
& \text{in }\mathbb{R}^1, \\
-\displaystyle\frac{1}{4\pi}\ln{(r^2+\epsilon^2)} +
  \displaystyle\frac{1}{4\pi}\sum_{\ell=1}^{n}
  \frac{1}{\ell}\epsilon^{2\ell}(r^2+\epsilon^2)^{-\ell}
& \text{in }\mathbb{R}^2, \\
\displaystyle
\frac{1}{4\pi}\sum_{\ell=0}^n {-\frac{1}{2}
    \choose \ell}\left(-\epsilon^2 \right )^\ell \left (r^2+\epsilon^2\right
  )^{-1/2-\ell}
& \text{in }\mathbb{R}^3,
\end{cases}
\end{align*}
where $\epsilon$ is the perturbation size.  The generalized binomial
coefficient is defined as
\begin{align*} 
  {\alpha \choose \ell} = 
  \frac{(\alpha)_{\ell}}{\ell!} = 
  \frac{1}{\ell!}\displaystyle\prod_{k=0}^{\ell-1}\left(\alpha-k\right),
\end{align*}
where $\alpha \in \mathbb{R}$ and $(\alpha)_{\ell}$ is the falling
factorial.  By construction, for any $\epsilon > 0$ and $r \neq 0$,
$G^{\epsilon,n}(r) \to G(r)$ pointwise as $n\to\infty$.  Note that
when $n=0$, the one-parameter family of regularized kernels,
equation~\eqref{eqn:regularized_n0}, is recovered.

\subsection{Hamiltonian of a regularized system}
\label{sec:error_particular}

If a time integrator is used to generate an approximate numerical
solution to the system
\begin{align}
  \ddot{\xx}_{j} = -\sum_{k \neq j} w_{jk}
    \nabla G^{\epsilon,n}(\xx_{j} - \xx_{k}),
  \label{eqn:regularized_system_two}
\end{align}
the resulting error in the Hamiltonian can be decomposed into two
parts: the modelling error that arises from replacing
$\nabla_{\mathbf{x}} G(\|\mathbf{x}-\mathbf{y}\|)$ with
$\nabla_{\mathbf{x}} G^{\epsilon,n}(\|\mathbf{x}-\mathbf{y}\|)$, and
the time-stepping error associated with discrete time integration.
The error in the Hamiltonian can be bounded,
\begin{align*}
  |H^{\epsilon,n}(t) - H(0)| &\leq |H^{\epsilon,n}(t) - H^{\epsilon,n}(0)| +
  |H^{\epsilon,n}(0) - H(0)|.
\end{align*}
The first term is the time-stepping error, $|H^{\epsilon,n}(t) -
H^{\epsilon,n}(0)|$, which depends on the chosen time step size.  For
$\epsilon>0$ and $n\ge0$, the time stepping error goes to zero as
$\Delta t \to 0$.  The rate at which this error term goes to zero
depends, however, on the size of the time step relative to $\nabla
G^{\epsilon,n}$.  The second term is the modelling or smoothing error,
and satisfies
\begin{align}
  |H^{\epsilon,n}(0) - H(0)| = \left| \sum_{j} \sum_{k\neq j} w_{jk}
  \left(G^{\epsilon,n}(\xx_{j} - \xx_{k}) - G(\xx_{j} - \xx_{k}) 
  \right) \right|.
  \label{eqn:modelling_error}
\end{align}
In the numerical experiments in Section~\ref{sec:numerics}, the
quantity $|H^{\epsilon,n}(T) - H(0)|$ will be reported for various
choices of $\epsilon$, $n$, and $\Delta t$.

\subsection{Global smoothing errors} Since the regularized system,
equation~\eqref{eqn:regularized_system_two}, relies on the gradient of
the regularized kernels, we investigate the error of the gradient.  One
way to quantify the quality of the regularization is to measure the {\em
global smoothing error},
\begin{align}
  e[\epsilon,n] = 
  \begin{cases}
    \displaystyle
    2\int_{0}^\infty \left|\frac{\partial
      G^{\epsilon,n}(r)}{\partial r} - \frac{\partial G(r)}{\partial
      r}\right|\,dr
    & \text{ in } \mathbb{R}^1, 
    \vspace*{0.2cm} \\
    \displaystyle
    \int_{0}^\infty \left|\frac{\partial
      G^{\epsilon,n}(r)}{\partial r} - \frac{\partial G(r)}{\partial
      r} \right|2\pi r\,dr 
    & \text{ in } \mathbb{R}^2,
    \vspace*{0.2cm} \\
    \displaystyle
    \int_{0}^\infty \left|\frac{\partial
      G^{\epsilon,n}(r)}{\partial r} - \frac{\partial G(r)}{\partial
      r} \right|4\pi r^2\,dr 
    & \text{ in } \mathbb{R}^3,
  \end{cases}
  \label{eqn:total_err}
\end{align}
where the gradient of the regularized kernels are
\begin{align}
  \frac{\partial G^{\epsilon,n}(r)}{\partial r} =
  \begin{cases}
    -\displaystyle r\sum_{\ell=0}^n {\frac{1}{2}
      \choose \ell}\left(-\epsilon^2 \right )^\ell
      \left(\frac{1}{2}-\ell\right)\left (r^2+\epsilon^2\right
    )^{-1/2-\ell}
    & \text{in }\mathbb{R}^1,\\
    \displaystyle
    -\frac{r}{2\pi}\sum_{\ell=0}^{n}
    \epsilon^{2\ell}(r^2+\epsilon^2)^{-\ell-1}
    & \text{in }\mathbb{R}^2,\\
    \displaystyle
    -\frac{r}{2\pi}\sum_{\ell=0}^n {-\frac{1}{2}
      \choose \ell}\left(-\epsilon^2 \right )^\ell  
    \left(-\frac{1}{2}-\ell\right)\left(r^2+\epsilon^2\right
    )^{-3/2-\ell}
    & \text{in }\mathbb{R}^3.
  \end{cases}
  \label{eqn:gradG_regularized}
\end{align}
If $r=0$, $\frac{\partial G^{\epsilon,n}(r)}{\partial r}=0$ which is
consistent with the physical argument that a particle does not feel any
self-force.   In $\mathbb{R}^{2}$, the expression for
$\frac{\partial}{\partial r}G^{\epsilon,n}(r)$ is a geometric series,
leading to the simplified expression
\begin{align}
  \frac{\partial G^{\epsilon,n}(r)}{\partial r} = -\frac{1}{2\pi r}
  \left(1 - \left(\frac{\epsilon^{2}}{\epsilon^{2}+r^2}\right)^{n+1}\right).
  \label{eqn:gradG_regularized_R2}
\end{align}
While this results in a more efficient expression for the gradient of
the regularized kernel, similar simplifications in $\mathbb{R}^{1}$ and
$\mathbb{R}^{3}$ do not exist.  The difference between the gradient of
the regularized and non-regularized kernels is 
\begin{align*}
 \frac{\partial}{\partial r} &\left(
 G^{\epsilon,n}(r) - G(r) \right) =
 \begin{cases}
   \displaystyle
   \sum_{\ell=n+1}^{\infty}
   2r\left(\frac 12 -\ell\right){\frac12 \choose
     \ell}(-1)^\ell\epsilon^{2i}
   (r^2+\epsilon^2)^{-\frac 12 - \ell}
   & \text{in }\mathbb{R}^1,\\
    \displaystyle
    \frac{1}{2\pi r}\left(\frac{\epsilon^{2}}{\epsilon^{2} + r^{2}}
    \right)^{n+1}
    & \text{in }\mathbb{R}^2,\\
    \displaystyle
    \frac{r}{2\pi}
    \sum_{\ell=n+1}^{\infty}
        {-\frac12 \choose \ell}\left(-\epsilon^2 \right)^\ell  
        \left(-\frac12-\ell\right)\left(r^2+\epsilon^2\right)^{-3/2-\ell}
        & \text{in }\mathbb{R}^3,
 \end{cases}
\end{align*}
where we have used equation~\eqref{eqn:gradG_regularized_R2} to
eliminate the summation in $\mathbb{R}^{2}$. The smoothing error,
equation~\eqref{eqn:total_err}, is computed by integrating these
expressions.  In Section~\ref{sec:analysis}, we provide estimates and
bounds for the global smoothing error.

\begin{figure}[htbp]
  \centering`<
  \subfloat[Pointwise error in
  $\mathbb{R}^1$]{\includegraphics[scale=1]{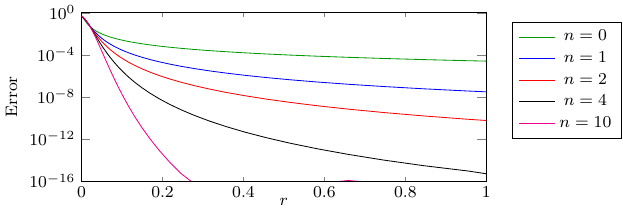}} \\
  \subfloat[Pointwise error in
  $\mathbb{R}^2$]{\includegraphics[scale=1]{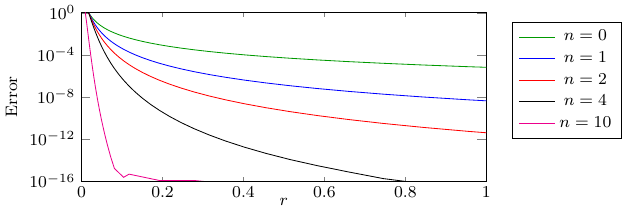}} \\ 
  \subfloat[Pointwise error in 
  $\mathbb{R}^3$]{\includegraphics[scale=1]{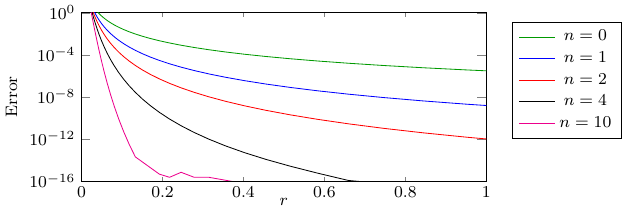}}
  \caption{Pointwise error of $\frac{\partial}{\partial r}
    (G^{\epsilon,n}-G)$ in $\mathbb{R}^1$, $\mathbb{R}^2$ and
    $\mathbb{R}^3$ as a function of $r$ for various $(\epsilon,n)$
    pairings.  The integral under each of these error curves is kept
    fixed.}
  \label{fig:constant_global_error}
\end{figure}

We are interested in pairings $(\epsilon,n)$ that give rise to kernels
with the same global smoothing error, equation~\eqref{eqn:total_err}.
The global smoothing error can be approximated by truncating the
infinite integral at 1, i.e.
\begin{align}
  \displaystyle\int_{B(0,1)} \left|
    \nabla G^{\epsilon,n} - \nabla G\right|\,d\mathbf{x} = 10^{-2},
  \label{eqn:global_smoothing_error}
\end{align}
where $B(0,1)$ is the unit ball centered at the origin in the
appropriate dimension.  This will facilitate a fair comparison (same
global smoothing error) between lower-order (small $n$) and higher-order
(large $n$) kernels.  Figure~\ref{fig:constant_global_error} shows the
pointwise error, equation~\eqref{eqn:gradG_regularized}, for various
kernels that satisfy equation~\eqref{eqn:global_smoothing_error}.  For
larger values of $n$, the error in the far field is greatly reduced
without introducing sharp derivatives near the singularity.

The pairings $(\epsilon,n)$ used to generate the curves with fixed
global smoothing error in Figure~\ref{fig:constant_global_error} are
summarized in Table~\ref{tbl:n_eps_pairs2}.  Instead of computing the
global smoothing error exactly, we satisfy
equation~\eqref{eqn:global_smoothing_error} by applying the trapezoid
rule on a sufficiently refined grid.  By choosing these $(\epsilon,n)$
pairings, the numerical experiments in Section~\ref{sec:numerics} will
demonstrate that although the global smoothing error (i.e.~a measure of
the error in the gradient) is held fixed, the modelling error,
equation~\eqref{eqn:modelling_error}, can be decreased with larger
values of $n$.
\begin{SCtable}
  \centering
  \begin{tabular}{cccc}
    $n$ & $\epsilon$ ($\mathbb{R}^{1})$ & $\epsilon$ ($\mathbb{R}^{2})$ 
    & $\epsilon$ ($\mathbb{R}^{3}$) \\
    \hline
    0  & $1.0051\times 10^{-2}$ & $6.3923\times 10^{-3}$ & 
         $5.0189\times 10^{-3}$ \\
    1  & $2.0001\times 10^{-2}$ & $1.2733\times 10^{-2}$ & 
         $1.0001\times 10^{-2}$ \\
    2  & $2.6667\times 10^{-2}$ & $1.6977\times 10^{-2}$ & 
         $1.3333\times 10^{-2}$ \\
    4  & $3.6572\times 10^{-2}$ & $2.3283\times 10^{-2}$ & 
         $1.8286\times 10^{-2}$ \\
    10 & $5.6755\times 10^{-2}$ & $3.6132\times 10^{-2}$ & 
         $2.8378\times 10^{-2}$
  \end{tabular}
  \caption{Values of $\epsilon$ and $n$ that result in a global
    smoothing error of $10^{-2}$.  These pairings will be used in
    Figures~\ref{fig:r1_oscillator_error}, \ref{fig:phase_plane},
    \ref{fig:r2_oscillator_error}, \ref{fig:r3_oscillator_error}, and
    \ref{fig:periodic_orbits}.}
  \label{tbl:n_eps_pairs2} 
\end{SCtable}

\section{Validity of the two-parameter family of regularized kernels}
\label{sec:validity}

This section demonstrates that the two-parameter family of regularized
kernels tend to the fundamental solutions for the Laplacian operator.
Specifically, we will show that for any $\epsilon > 0$, $\Delta
G^{\epsilon,n}(\xx) \to \Delta G(\xx) = -\delta(\xx)$ as $n\to\infty$,
and that for any $n$, $\Delta G^{\epsilon,n}(\xx) \to -\delta(\xx)$ as
$\epsilon\to0$, where $\delta$ is the delta function. All the
convergences are understood to be taken in the weak sense.  To check
this convergence, we will show that $\Delta G^{\epsilon,n}(\xx) \to 0$
for all $\xx \neq 0$, and that $\int \Delta G^{\epsilon,n}(\xx)d\xx$
converges to $-1$.  This implies that $G^{\epsilon,n}(r)$ converges to
the fundamental solution of the Laplacian operator.

\begin{remark}
  In general, if $G^{\epsilon,n}(\xx)$ converges pointwise to $G(\xx)$
  for all $\xx \neq 0$, then $\Delta G^{\epsilon,n}$ converges weakly
  to the delta function.  However, we will prove this result directly
  for our particular choice of $G^{\epsilon,n}$ as this will lead to
  closed-form solutions of $\Delta G^{\epsilon,n}$ which may be useful
  in future work.
\end{remark}

We first generate simplified expressions for $\Delta G^{\epsilon,n}$,
which will simplify the proofs that the two-parameter family of regularized kernels tend to
the fundamental solutions of the Laplace operator.

\subsection{Simplified expressions for $\Delta G^{\epsilon,n}(r)$}
Applying the Laplacian operator to the two-parameter family of
regularized kernels gives
\begin{align*}
  \Delta G^{\epsilon,n}(r) = 
  \begin{cases}
    \displaystyle
    -\sum_{\ell=0}^n {\frac{1}{2} \choose \ell}
    (-\epsilon^2)^\ell\left(\frac{1}{2}-\ell\right)
    (r^2+\epsilon^2)^{-\frac{3}{2}-\ell}(-2\ell r^2+\epsilon^2)
    & \text{in }\mathbb{R}^1, \\
    \displaystyle
    -\frac{1}{\pi}
    \sum_{\ell=0}^{n}\epsilon^{2\ell}(r^2+\epsilon^2)^{-\ell-2}(-\ell r^2+\epsilon^2)
    & \text{in }\mathbb{R}^2, \\
    \displaystyle
    -\frac{1}{2\pi}\sum_{\ell=0}^n {-\frac{1}{2} \choose \ell}
    (-\epsilon^2)^\ell\left(\frac{1}{2}+\ell\right)
    (r^2+\epsilon^2)^{-\frac{5}{2}-\ell}(-2\ell r^2+3\epsilon^2)
    & \text{in }\mathbb{R}^3.
  \end{cases}
\end{align*}
In $\mathbb{R}^1$, the generalized binomial coefficients can be
eliminated by using the identity
\begin{align*}
   {\frac{1}{2} \choose \ell}\left(\frac{1}{2}-\ell\right) = 
   \frac{(-1)^\ell\Gamma\left(\ell+\frac{1}{2}\right)}{2\sqrt{\pi}\Gamma(\ell+1)}.
\end{align*}
Thus in $\mathbb{R}^1$, the Laplacian operator applied to the
two-parameter family of regularized kernels is
\begin{align*}
  \Delta G^{\epsilon,n}(r) &= 
  -\frac{1}{2\sqrt{\pi}}\displaystyle
  \sum_{\ell=0}^n 
  \frac{\Gamma\left(\ell+\frac{1}{2}\right)}{\Gamma(\ell+1)}
  \epsilon^{2\ell}
  (r^2+\epsilon^2)^{-\frac{3}{2}-\ell}(-2\ell r^2+\epsilon^2)\\
  &=-\frac{1}{2\epsilon\sqrt{\pi}}\displaystyle
  \sum_{\ell=0}^n 
  \frac{\Gamma\left(\ell+\frac{1}{2}\right)}{\Gamma(\ell+1)}
  \frac{  \left(-2\ell
  a^2+1\right)}{\left(a^2+1\right)^{\frac{3}{2}+\ell}}
\end{align*}
where $a = \frac{r}{\epsilon}$.  Similarly, using the identity
\begin{align*}
   {-\frac{1}{2} \choose \ell}\left(\frac{1}{2}+\ell\right) = 
   \frac{(-1)^\ell\Gamma\left(\ell+\frac{3}{2}\right)}{\sqrt{\pi}\Gamma(\ell+1)},
\end{align*}
the Laplacian of the two-parameter family of regularized kernels in $\mathbb{R}^3$ is
\begin{align*}
 \Delta G^{\epsilon,n}(r) = 
  -\frac{1}{2\pi\sqrt{\pi}\epsilon^3}\displaystyle
  \sum_{\ell=0}^n 
  \frac{\Gamma\left(\ell+\frac{3}{2}\right)}{\Gamma(\ell+1)}
  \frac{  \left(-2\ell
  a^2+3\right)}{\left(a^2+1\right)^{\frac{5}{2}+\ell}}.
\end{align*}
These expressions for $\Delta G^{\epsilon,n}$in
$\mathbb{R}^{1}$, $\mathbb{R}^{2}$, and $\mathbb{R}^{3}$ can be
simplified further so that there is no summation.  These much simpler
expressions allow us to verify that $\Delta G^{\epsilon,n}$ converges
weakly to the delta function.
\begin{theorem}
The Laplacian operator applied to the two-parameter family of
regularized kernels can be expressed as 
  \begin{align}
    \Delta G^{\epsilon,n}(r) = 
    \begin{cases}
      \displaystyle
      -\frac{1}{\epsilon\sqrt{\pi}}
      \left(\frac{1}{1+ a^2}\right)^{n+\frac{3}{2}}
      \frac{\Gamma\left(n+\frac{3}{2}\right)}{\Gamma(n+1)}
      & \text{in }\mathbb{R}^1, \\
      \displaystyle
      -\frac{(n+1)}{\pi\epsilon^2}
      \left(\frac{1}{1+a^2}\right)^{n+2}
      & \text{in }\mathbb{R}^2, \\
      \displaystyle
      -\frac{1}{\epsilon^3\pi\sqrt{\pi}}
      \left(\frac{1}{1+ a^2}\right)^{n+\frac{5}{2}}
      \frac{\Gamma\left(n+\frac{5}{2}\right)}{\Gamma(n+1)}
      & \text{in }\mathbb{R}^3,
    \end{cases}
    \label{eqn:gxx}
  \end{align}
where $a = \frac{r}{\epsilon}$.
\end{theorem}
\begin{proof}
  In $\mathbb{R}^1$, since $\Gamma(z+1) = z\Gamma(z)$,
  equation~\eqref{eqn:gxx} holds for $n=0$.  Now, suppose
  equation~\eqref{eqn:gxx} holds for $n = p$.  Then,
  \begin{align*}
    \Delta G^{\epsilon,p+1}(r) &= -\frac{1}{\epsilon\sqrt{\pi}}
    \left(\frac{1}{1+ a^2}\right)^{p+\frac{3}{2}}
    \frac{\Gamma\left(p+\frac{3}{2}\right)}{\Gamma(p+1)}
    -
    \frac{\Gamma\left(p+\frac{3}{2}\right)}{2\epsilon\sqrt{\pi}\Gamma(p+2)}
    \frac{\left(-2(p+1)a^2+1\right)}
    {\left(a^2+1\right)^{\frac{5}{2}+p}} \\
    &= -\frac{1}{\epsilon\sqrt{\pi}}
    \left(\frac{1}{1+ a^2}\right)^{p+\frac{5}{2}}
    \frac{\Gamma\left(p+\frac{3}{2}\right)}{\Gamma(p+2)}
    \left[
      (1+a^2)\frac{\Gamma(p+2)}{\Gamma(p+1)}
      + \frac{\left(-2(p+1)a^2+1\right)}{2}
    \right] \\
    &= -\frac{1}{\epsilon\sqrt{\pi}}
    \left(\frac{1}{1+ a^2}\right)^{p+\frac{5}{2}}
    \frac{\Gamma\left(p+\frac{3}{2}\right)}{\Gamma(p+2)}
    \left[
      (1+a^2)(p+1)
      + \frac{\left(-2(p+1)a^2+1\right)}{2}
    \right]\\
    & = -\frac{1}{\epsilon\sqrt{\pi}}
    \left(\frac{1}{1+ a^2}\right)^{p+\frac{5}{2}}
    \frac{\Gamma\left(p+\frac{3}{2}\right)}{\Gamma(p+2)}
    \left(p+\frac{3}{2}\right)\\
    &= -\frac{1}{\epsilon\sqrt{\pi}}
    \left(\frac{1}{1+ a^2}\right)^{p+\frac{5}{2}}
    \frac{\Gamma\left(p+\frac{5}{2}\right)}{\Gamma(p+2)},
  \end{align*}
  which establishes that the equivalence must hold for $n=p+1$.  A
  similar inductive argument can be used to establish
  equation~\eqref{eqn:gxx} in $\mathbb{R}^2$ and $\mathbb{R}^3$.
  \qed
\end{proof}
\begin{remark}
  Since closed-form expressions (without the infinite summations) can
  be obtained for $\Delta G^{\epsilon,n}(r)$, one might expect that
  similar closed form expressions should exist for $\frac{\partial
  G^{\epsilon,n}(r)}{\partial r}$.  This is unfortunately not the case
  in $\mathbb{R}^1$ and $\mathbb{R}^3$, as previously noted.
\end{remark}



\subsection{Weak convergence to the delta function}
Now that a closed-form expression for the Laplacian of the
two-parameter family of regularized kernels has been attained, their
weak limits can be shown to converge to the delta function.
\begin{theorem}
  The Laplacian of the two-parameter family of regularized kernels,
  $G^{\epsilon,n}(r)$, converges weakly to the delta function.  That is,
  \begin{align}
    \label{eqn:lim1}
    \lim_{n\to\infty}\Delta G^{\epsilon,n}(r) = -\delta(r) \text{ for any }
    \epsilon  > 0, \\
    \label{eqn:lim2}
    \lim_{\epsilon\to0}\Delta G^{\epsilon,n}(r) = -\delta(r) \text{ for any } n \ge 0.
  \end{align}
\end{theorem}
\begin{proof}
  We first check that the two limits, equation~\eqref{eqn:lim1} and
  equation~\eqref{eqn:lim2} converge to 0 for all $r \neq 0$.  If
  $\epsilon>0$ and $r \neq 0$, then, from equation~\eqref{eqn:gxx},
  \begin{align*}
    \lim_{n\to\infty}\Delta G^{\epsilon,n}(r) = 0,
  \end{align*}
  for $\mathbb{R}^1, \mathbb{R}^2$ and $\mathbb{R}^3$.  Next, since
  equation~\eqref{eqn:gxx} can be rewritten as
  \begin{align*}
    \Delta G^{\epsilon,n}(r) = 
    \begin{cases}
      \displaystyle
      -\frac{\epsilon^{2n+2}}{\sqrt{\pi}}
      \left(\frac{1}{\epsilon^2+ r^2}\right)^{n+\frac{3}{2}}
      \frac{\Gamma\left(n+\frac{3}{2}\right)}{\Gamma(n+1)}
      & \text{in }\mathbb{R}^1,\\
      \displaystyle
      -\frac{(n+1)\epsilon^{2n+2}}{\pi}
      \left(\frac{1}{\epsilon^2+r^2}\right)^{n+2}
      & \text{in }\mathbb{R}^2,\\
      \displaystyle
      -\frac{\epsilon^{2n+2}}{\pi\sqrt{\pi}}
      \left(\frac{1}{\epsilon^2+ r^2}\right)^{n+\frac{5}{2}}
      \frac{\Gamma\left(n+\frac{5}{2}\right)}{\Gamma(n+1)}
      & \text{in }\mathbb{R}^3,
    \end{cases}
  \end{align*}
  it follows that for any $n \ge 0$, $r \neq 0$,
  \begin{align*}
    \lim_{\epsilon\to0}\Delta G^{\epsilon,n}(r) = 0.
  \end{align*}
  Next, for all $\epsilon$ and $n$,
  \begin{align*}
    2\int_0^\infty-\frac{1}{\epsilon\sqrt{\pi}}\left(\frac{1}{1+\frac{r^2}{\epsilon^2}}\right)^{n+\frac{3}{2}}
    \frac{\Gamma(n+\frac{3}{2})}{\Gamma(n+1)}\,dr = -1, & \quad \text{in
    }\mathbb{R}^1, \\
    -\frac{n+1}{\pi\epsilon^2}\int_0^\infty
    \left(\frac{1}{1+\frac{r^2}{\epsilon^2}}\right)^{n+2}2\pi r\,dr =
    -1, & \quad \text{in }\mathbb{R}^2, \\
    -\frac{1}{\epsilon^3\pi\sqrt{\pi}}\int_0^\infty\left(\frac{1}{1+\frac{r^2}{\epsilon^2}}\right)^{n+\frac{5}{2}}
    \frac{\Gamma(n+\frac{5}{2})}{\Gamma(n+1)}4\pi r^2\,dr = -1, & \quad
    \text{in }\mathbb{R}^3. 
  \end{align*}
  Finally, let $f$ be a compactly supported smooth function and
  $\epsilon > 0$.   Then, for any $\alpha > 0$,
  \begin{align*}
    \int_{\mathbb{R}^{d}}\lim_{n \to \infty}
      \Delta G^{\epsilon,n}(\xx)f(\xx)d\xx = 
    \int_{B(0,\alpha)}\lim_{n \to \infty}
      \Delta G^{\epsilon,n}(\xx)f(\xx)d\xx.
  \end{align*}
  Since this holds for all $\alpha > 0$, we have
  \begin{align*}
    \int_{\mathbb{R}^{d}}\lim_{n \to \infty} 
      \Delta G^{\epsilon,n}(\xx)f(\xx)d\xx &= 
    f(0)\int_{B(0,\alpha)} \lim_{n \to \infty}
      \Delta G^{\epsilon,n}(\xx)d\xx \\
    &=f(0)\int_{\mathbb{R}^{d}} \lim_{n \to \infty}
      \Delta G^{\epsilon,n}(\xx)d\xx \\
    &=-f(0),
  \end{align*}
  which establishes equation~\eqref{eqn:lim1}.  The proof of
  equation~\eqref{eqn:lim2} is similarly proved by fixing $n>0$ and
  showing that
  \begin{align*}
    \int_{\mathbb{R}^{d}}\lim_{\epsilon \to 0}
      \Delta G^{\epsilon,n}(\xx)f(\xx)d\xx = -f(0).
  \end{align*}
  \qed
\end{proof}

\section{Error analysis}
\label{sec:analysis}
In this section we provide estimates of the global smoothing error,
equation~\eqref{eqn:total_err}, that arises from using the
two-parameter family of regularized kernels in $\mathbb{R}^1,
\mathbb{R}^2$ and $\mathbb{R}^3$.  In all three dimensions, we assume
that the point masses are contained in the ball of radius $R$ centered
at the origin.

\subsection{Global smoothing error in $\mathbb{R}^1$}
We start by bounding 
\begin{align}
  S[\epsilon,n] = 2\sum_{\ell=n+1}^{\infty}
  {\frac12 \choose \ell}(-1)^\ell \left(
  z^{\ell-\frac12} -1 \right),
  \label{eqn:r1_sne}
\end{align}
where $z = \frac{\epsilon^{2}}{R^{2}+\epsilon^{2}} \in (0,1)$, which
will arise in our error estimate shortly.

\begin{theorem}
  \label{thm:r1_S}
  Let $S$ be defined in equation~\eqref{eqn:r1_sne}.  For any
  $\epsilon > 0$,
  \begin{align*}
    \lim_{n \to \infty} S[\epsilon,n] = 0.
  \end{align*}
\end{theorem}
\begin{proof}
  Since $z \in (0,1)$, we can bound $S[\epsilon,0]$ as
  \begin{align*}
    S[\epsilon,0] &=2 \sum_{\ell=1}^{\infty}
    {\frac12 \choose \ell}(-1)^\ell\left(
      z^{\ell-\frac12}-1 \right) \\
    &= 2 \left(\frac{-\sqrt{z}}{1+\sqrt{1-z}} + 1\right) \in (0,2).
  \end{align*}
  Therefore, since $S[\epsilon,0]$ is bounded, this gives the desired result, that
  \begin{align*}
    \lim_{n \to \infty} S[\epsilon,n] = 0.
  \end{align*}
  \qed
\end{proof}
Next, the value of $S[\epsilon,n]$ is estimated.  Define
\begin{align*}
  g(z) := z^{-\frac12}\left(
    S[\epsilon,0] - 2\right) =  \frac{-2}{1 + \sqrt{1-z}}
       = 2\sum_{\ell=1}^{\infty} {\frac12 \choose \ell}  
       (-1)^\ell z^{\ell-1}.
\end{align*}
By Taylor's theorem, we have
\begin{align}
  2\sum_{\ell=n+1}^{\infty}{\frac12 \choose \ell}(-1)^\ell
  z^{\ell-1} = 
  \frac{-2}{1+\sqrt{1-z}} - 
  2\sum_{\ell=1}^{n}    {\frac12 \choose \ell}(-1)^\ell z^{\ell-1} 
  = \frac{g^{(n)}(\xi)}{n!}z^n,
  \label{eqn:g_r1}
\end{align}
where $\xi \in [0,z] \subset [0,1)$.  Multiplying both sides of
  equation~\eqref{eqn:g_r1} by $z^{\frac12}$ gives,
\begin{align}
  2\sum_{\ell=n+1}^{\infty}    
  {\frac12 \choose \ell}(-1)^\ell z^{\ell-\frac12} = 
  \frac{g^{(n)}(\xi)}{n!}z^{n+\frac12}.
  \label{eqn:r1_sne_term1}
\end{align}
Next, since
\begin{align*}
  2\sum_{\ell=1}^{\infty} {\frac12 \choose \ell}(-1)^\ell = -2,
\end{align*}
this gives
\begin{align}
  2\sum_{\ell=n+1}^{\infty} {\frac12 \choose \ell}(-1)^\ell = -2 -
  2\sum_{\ell=1}^{n} {\frac12 \choose \ell}(-1)^\ell .
  \label{eqn:r1_sne_term2}
\end{align}
Substituting equation~\eqref{eqn:r1_sne_term1} and
equation~\eqref{eqn:r1_sne_term2} into equation~\eqref{eqn:r1_sne}
gives 
\begin{align*}
  S[\epsilon,n] = \left( 
    \frac{g^{(n)}(\xi)}{n!}z^{n+\frac12} + 2\left(
    1 + \sum_{\ell=1}^{n} {\frac12 \choose \ell}(-1)^\ell 
    \right)\right).
\end{align*}
The global smoothing error~\eqref{eqn:total_err} can be now be estimated
as
\begin{align*}
  e[\epsilon,n] &= 2 \int_0^R \left|\frac{\partial}{\partial r} \left(
  G^{\epsilon,n}(r) -G(r) \right)\right|\,dr \\
  &= 2 \int_0^R \frac{\partial}{\partial r} \left(
  \sum_{\ell=n+1}^{\infty}
 {\frac12 \choose \ell}(-1)^\ell\epsilon^{2\ell}
 (r^2+\epsilon^2)^{\frac 12 - \ell} \right) \,dr \\
  &=  2 \left.\left(
  \sum_{\ell=n+1}^{\infty}
 {\frac12 \choose \ell}(-1)^\ell\epsilon^{2\ell}
 (r^2+\epsilon^2)^{\frac 12 - \ell} \right) \right|_0^R \\
  &= 2\epsilon\sum_{\ell=n+1}^{\infty}
 {\frac12 \choose \ell}(-1)^\ell
 \left(\left(\frac{\epsilon^2}{R^2+\epsilon^2}
 \right)^{\ell-\frac12} -1 \right) \\
  &= \epsilon S[\epsilon,n] \\
  &= \epsilon \left(
    \frac{g^{(n)}(\xi)}{n!}
    \left(\frac{\epsilon^2}{R^2+\epsilon^2}\right)^{n+\frac12}
    + 2\left(1 + \sum_{\ell=1}^{n} {\frac12 \choose \ell}(-1)^\ell 
    \right)\right),
\end{align*}
where the absolute value was dropped since the derivative of each
term in the summation is positive.  This provides us with an estimate
of the global smoothing error, and as expected, Theorem~\ref{thm:r1_S}
guarantees that
\begin{align*}
  \text{If } \epsilon>0, \lim_{n\to\infty}e[\epsilon,n] &= 0, \\
  \text{if } n \geq 0, \lim_{\epsilon\to 0}e[\epsilon,n] &= 0.
\end{align*}


\subsection{Global smoothing error in $\mathbb{R}^2$}
Following the analysis in $\mathbb{R}^1$, we start by bounding
\begin{align}
  \nonumber
  S[\epsilon,n] &= \sum_{\ell=n+1}^{\infty} \left(-\frac{1}{2\ell-1}
    \left(\frac{z}{1-z}\right)^{\ell-\frac12}
    \tensor[_2]{F}{_1}\left(\ell+1,\ell-\frac12;\ell+\frac12;
    -\frac{z}{1-z}\right) \right. \\
    \label{eqn:S_r2}
    &\qquad\qquad\qquad\qquad-\left.
    \frac{(-1)^{\ell}\pi^{\frac32}}{4\ell!\Gamma(\frac32-\ell)}\right),
\end{align}
where $z = \frac{\epsilon^{2}}{R^{2}+\epsilon^{2}} \in (0,1)$, which
will arise in the error estimate for $\mathbb{R}^2$ .  The
hypergeometric function is defined as
\begin{align*}
  \tensor[_2]{F}{_1}(a,b;c;z) =
  \sum_{n=0}^{\infty}\frac{(a)_{n}(b)_{n}}{(c)_{n}}\frac{z^{n}}{n!},
\end{align*}
where $(a)_{n} = a(a-1)(a-2)\cdot(a-n+1)$ is the falling factorial, and
$|z|<1$.
\begin{theorem}
  \label{thm:r2_S}
  Let $S$ be defined in equation~\eqref{eqn:S_r2}.  For any $\epsilon > 0$, 
  \begin{align*}
    \lim_{n \to \infty} S[\epsilon,n] = 0.
  \end{align*}
\end{theorem}
\begin{proof}
Unlike in $\mathbb{R}^{1}$, we can not eliminate the summation in
$S[\epsilon,n]$, but it can be partially simplified as
\begin{align}
  S[\epsilon,0] = \frac{\pi}{2}-\sum_{\ell=1}^{\infty}\frac{1}{2\ell-1}
  \left(\frac{z}{1-z}\right)^{\ell-\frac12} \tensor[_2]{F}{_1}
  \left(\ell+1,\ell-\frac12;\ell+\frac12,-\frac{z}{1-z}\right).
  \label{eqn:r2_S}
\end{align}
In Figure~\ref{fig:2d_bound}, we plot the 150 term partial sum of
equation~\eqref{eqn:r2_S} which guarantees five digits of accuracy.
Therefore, since $S[\epsilon,0] \in (0,\frac{\pi}{2})$ is bounded, this
gives the desired result, that
\begin{align*}
  \lim_{n \to \infty} S[\epsilon,n] = 0.
\end{align*}
\begin{figure}[htbp]
  \begin{minipage}{0.7\textwidth}
    \centering
    \includegraphics[scale=1]{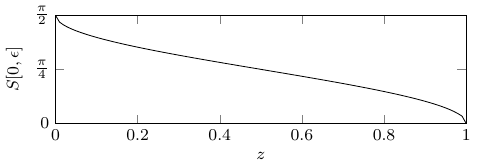}
  \end{minipage}\hfill
  \begin{minipage}{0.29\textwidth}
    \caption{\label{fig:2d_bound} A partial sum of
    equation~\eqref{eqn:r2_S}.  As desired,  $S[\epsilon,0]$ is bounded
    for $z \in (0,1).$}
  \end{minipage}
\end{figure}
\qed
\end{proof}
Without a closed-form expression for $S[\epsilon,n]$, we are unable to
use Taylor's theorem to provide an estimate of $S[\epsilon,n]$.
However, values of $S[\epsilon,n]$ can be computed by taking partial
sums of equation~\eqref{eqn:r2_S}.  The global smoothing error,
equation~\eqref{eqn:total_err}, can be expressed in terms of
$S[\epsilon,n]$\footnote{If the simplified expression for
$\frac{\partial G^{\epsilon,n}}{\partial r}$ in
equation~\eqref{eqn:gradG_regularized_R2} is used, one recovers the
same expression involving the hypergeomtric function.},
\begin{align*}
  e[\epsilon,n] &= 2\pi \int_{0}^{R}r\left| 
    \frac{\partial}{\partial r} \left(G^{\epsilon,n}(r) - G(r)\right)
  \right|\,dr \\
  &= \int_{0}^{R} r^{2} \left| \sum_{\ell=n+1}^{\infty} 
    \epsilon^{2\ell}(r^{2}+\epsilon^2)^{-\ell-1} \right|\,dr \\
  &=\epsilon \sum_{\ell=n+1}^{\infty}\left\{-\frac{1}{2\ell-1}
    \left(\frac{\epsilon}{R}\right)^{2\ell-1}
    \tensor[_2]{F}{_1}\left(\ell+1,\ell-\frac12;
    \ell+\frac12,-\frac{\epsilon^{2}}{R^2}\right)\right. \\
  &\qquad\qquad\qquad\qquad\left.
    -\frac{(-1)^{\ell}\pi^{\frac32}}{4\ell!\Gamma(\frac32-\ell)}
  \right\} \\
  &=\epsilon \sum_{\ell=n+1}^{\infty}\left\{-\frac{1}{2\ell-1}
    \left(\frac{z}{1-z}\right)^{\ell-\frac12}
    \tensor[_2]{F}{_1}\left(\ell+1,\ell-\frac12;
    \ell+\frac12,-\frac{z}{1-z}\right)\right. \\
  &\qquad\qquad\qquad\qquad\left.
    -\frac{(-1)^{\ell}\pi^{\frac32}}{4\ell!\Gamma(\frac32-\ell)}
  \right\} \\
  &=\epsilon S[\epsilon,n],
\end{align*}
As expected, Theorem~\ref{thm:r2_S} guarantees that
\begin{align*}
  \text{If } \epsilon>0,\lim_{n\to\infty}e[\epsilon,n] &= 0, \\
  \text{If } n \geq 0, \lim_{\epsilon\to 0}e[\epsilon,n] &= 0.
\end{align*}

\subsection{Global smoothing error in $\mathbb{R}^3$}
We start by bounding
\begin{align}
  S[\epsilon,n] &= -\frac{2n+3}{2(n+1)}{-\frac12 \choose n}(-1)^{n}
    \left(z^{n+\frac12}-1\right) \nonumber \\
    &\qquad\qquad -\frac12 \sum_{\ell=n+1}^{\infty}
    {-\frac12 \choose \ell} \frac{(-1)^{\ell}}{\ell+1}\left(
    z^{\ell+\frac12} - 1\right),
  \label{eqn:r3_sne}
\end{align}
where $z = \frac{\epsilon^{2}}{R^{2}+\epsilon^{2}} \in (0,1)$, which
will arise in the error estimate.

\begin{theorem}
  \label{thm:r3_S}
  Let $S$ be defined in equation~\eqref{eqn:r3_sne}.  For any
  $\epsilon > 0$,
  \begin{align*}
    \lim_{n \to \infty} S[\epsilon,n] = 0.
  \end{align*}
\end{theorem}
\begin{proof}
  Since $z \in (0,1)$, we can bound $S[\epsilon,0]$ as
  \begin{align*}
    S[\epsilon,0] &=-\frac{3}{2}\left(z^\frac12 -1 \right) - 
    \frac12\sum_{\ell=1}^{\infty} {-\frac12 \choose \ell}
    \frac{(-1)^{\ell}}{\ell+1}\left(z^{\ell+\frac12}-1 \right) \\
    &= 2 - \frac{3}{2}z^{\frac12} - \frac{z^{\frac{3}{2}}}
      {2\left(1+\sqrt{1-z}\right)^{2}} \in (0,2).
  \end{align*}
  Therefore, since $S[\epsilon,0]$ is bounded, this gives 
  \begin{align*}
    \lim_{n \to \infty} S[\epsilon,n] = 0.
  \end{align*}
  \qed
\end{proof}
To estimate the value of $S[\epsilon,n]$, define
\begin{align*}
  g(z) :=  -2z^{-\frac{3}{2}}\left(S[\epsilon,0]-2+
    \frac{3}{2}z^{\frac12} \right) = 
    \frac{1}{\left(1+\sqrt{1-z}\right)^{2}} = 
    \sum_{\ell=1}^{\infty} {-\frac12 \choose \ell}
    \frac{(-1)^{\ell}}{\ell+1}z^{\ell-1}.
\end{align*}
By Taylor's theorem, we have
\begin{align*}
  \sum_{\ell=n+1}^{\infty} {-\frac12 \choose \ell}
    \frac{(-1)^{\ell}}{\ell+1}z^{\ell-1} = 
    \frac{1}{\left(1+\sqrt{1-z}\right)^{2}} - 
    \sum_{\ell=1}^{n} {-\frac12 \choose \ell}
    \frac{(-1)^{\ell}}{\ell+1}z^{\ell-1}
  = \frac{g^{(n)}(\xi)}{n!}z^n,
\end{align*}
where $\xi \in [0,z] \subset [0,1)$.  Multiplying both sides of the
equation by $z^{\frac{3}{2}}$ gives,
\begin{align}
  \sum_{\ell=n+1}^{\infty} {-\frac12 \choose \ell}
    \frac{(-1)^{\ell}}{\ell+1}z^{\ell+\frac12} = 
    \frac{g^{(n)}(\xi)}{n!}z^{n+\frac{3}{2}}.
  \label{eqn:r3_sne_term1}
\end{align}
Next, since
\begin{align*}
  \sum_{\ell=1}^{\infty}    {-\frac12 \choose
  \ell}\frac{(-1)^\ell}{\ell+1} = 1,
\end{align*}
we have
\begin{align}
  \sum_{\ell=n+1}^{\infty}{-\frac12 \choose
  \ell}\frac{(-1)^\ell}{\ell+1} = 1 -
  \sum_{\ell=1}^{n} {-\frac12 \choose \ell}\frac{(-1)^\ell}{\ell+1}.
  \label{eqn:r3_sne_term2}
\end{align}
Substituting equation~\eqref{eqn:r3_sne_term1} and
equation~\eqref{eqn:r3_sne_term2} into equation~\eqref{eqn:r3_sne}, we
have
\begin{align*}
  S[\epsilon,n] &=  
  -\frac{2n+3}{2(n+1)}{-\frac12 \choose n}(-1)^{n} 
  \left(z^{n+\frac12}-1\right) -  \\
  & \qquad \frac12\frac{g^{(n)}(\xi)}{n!}z^{n+\frac32} + \frac{1}{2}
  \left(1 - \sum_{\ell=1}^{n}    {-\frac12 \choose
  \ell}(-1)^\ell\frac{1}{\ell+1} \right).
\end{align*}
The global smoothing error can be now be expressed in terms of
$S[\epsilon,n]$,
\begin{align*}
  e[\epsilon,n] &= 4\pi \int_0^R r^{2} \left|
    \frac{\partial}{\partial r} \left(
    G^{\epsilon,n}(r) -G(r) \right)\right|\,dr \\
  &= 4\pi \int_0^R r^{2} \left| \sum_{\ell=n+1}^{\infty}
    \frac{r}{2\pi}
    {-\frac12 \choose \ell}(-1)^{\ell}\epsilon^{2\ell}
    \left(-\frac12 -\ell\right)
    (r^{2} + \epsilon^{2})^{-\frac{3}{2}-\ell}
    \right|\,dr \\
  &= 2\sum_{\ell=n+1}^{\infty}{-\frac12 \choose \ell}(-1)^{\ell}
     \left(\frac12 + \ell\right) \epsilon^{2\ell} \int_{0}^{R}
     r^{3}(r^{2} + \epsilon^{2})^{-\frac{3}{2}-\ell}\,dr \\
  &= \sum_{\ell=n+1}^{\infty}{-\frac12 \choose \ell}(-1)^{\ell}
     \left(\frac12 + \ell\right) \epsilon^{2\ell} 
     \int_{\epsilon^{2}}^{R^{2}+\epsilon^2}
     \left(u-\epsilon^{2}\right)u^{-\frac{3}{2}-\ell}du \\
  &= \epsilon\sum_{\ell=n+1}^{\infty}{-\frac12 \choose \ell}(-1)^{\ell}
     \left\{\frac{\frac12+\ell}{\frac12-\ell}
     \left(\left(
     \frac{\epsilon^{2}}{R^2+\epsilon^2}\right)^{\ell-\frac12}
     -1\right) \right.\\
     & \qquad\qquad\qquad\qquad\qquad\qquad
     +\left.\left(\left(
      \frac{\epsilon^{2}}{R^2+\epsilon^2}\right)^{\ell+\frac12}
      -1\right)
     \right\},
\end{align*}
where the absolute value can be dropped since each term in the
summation is positive.  By shifting indices and using standard
properties of the generalized binomial coefficient, we have
\begin{align*}
  e[\epsilon,n] &= \epsilon\sum_{\ell=n+1}^{\infty}
    {-\frac12 \choose \ell}(-1)^{\ell}
    \left(\left(
    \frac{\epsilon^{2}}{R^{2}+\epsilon^2}
    \right)^{\ell+\frac12}-1\right) \\
    &-\epsilon\sum_{\ell=n}^{\infty}
    \frac{2\ell+3}{2(\ell+1)}{-\frac12 \choose \ell}(-1)^{\ell}
    \left(\left(
    \frac{\epsilon^{2}}{R^{2}+\epsilon^{2}}
    \right)^{\ell+\frac12}-1\right) \\
  &=\epsilon\left\{-\frac{2n+3}{2(n+1)}{-\frac12 \choose n}(-1)^{n}
    \left(\left(\frac{\epsilon^{2}}{R^{2}+\epsilon^{2}}\right)^
    {n+\frac12}-1\right)\right. - \\
    &\left.-\frac{1}{2}\sum_{\ell=n+1}^{\infty}{-\frac12 \choose \ell}
    (-1)^{\ell}\frac{1}{\ell+1}\left(\left(
    \frac{\epsilon^{2}}{R^{2}+\epsilon^{2}}\right)^
    {\ell+\frac12} - 1\right)\right\} \\
    &=\epsilon S[\epsilon,n] \\
    &=\epsilon\left\{
  -\frac{2n+3}{2(n+1)}{-\frac12 \choose n}(-1)^{n} 
  \left(\left(\frac{\epsilon^{2}}{R^{2}+\epsilon^{2}}\right)^
  {n+\frac12}-1\right)\right. \\
  &\left.
  -\frac12\frac{g^{(n)}(\xi)}{n!}\left(
  \frac{\epsilon^{2}}{R^{2}+\epsilon^{2}}\right)^
  {n+\frac32} + \frac{1}{2}
  \left(1 - \sum_{\ell=1}^{n}    {-\frac12 \choose
  \ell}(-1)^\ell\frac{1}{\ell+1} \right)\right\}.
\end{align*}
This provides us with an estimate of the global smoothing
error.  Invoking Theorem~\ref{thm:r3_S}
guarantees that
\begin{align*}
  \text{If } \epsilon >0, \lim_{n\to\infty}e[\epsilon,n] &= 0,  \\
  \text{If } n \geq 0, \lim_{\epsilon\to 0}e[\epsilon,n] &= 0.
\end{align*}

\section{Numerical Examples}
\label{sec:numerics}

We first demonstrate the described behaviors in $\mathbb{R}^1,
\mathbb{R}^2$ and $\mathbb{R}^3$ for a dynamical system with two
particles of equal mass---one with a positive unit charge and the
other with a unit negative charge.  Two different sets of
$(\epsilon,n)$ pairings are used.  Those in
Table~\ref{tbl:n_eps_pairs2} keep the global smoothing error fixed at
$10^{-2}$, as described in
equation~\eqref{eqn:global_smoothing_error}.  The pairings reported in
Table~\ref{tbl:n_eps_pairs3} keep the modelling error,
equation~\eqref{eqn:modelling_error}, fixed at $4.89 \times 10^{-6}$.
Note that the modelling error depends on the initial condition whereas
the global smoothing error does not.  We conclude the numerical
experiments by simulating 5 particles with an initial condition that
results in a periodic orbit.

\begin{SCtable}
  \centering
  \begin{tabular}{cccc}
    $n$ & $\epsilon$ ($\mathbb{R}^{1})$ & $\epsilon$ ($\mathbb{R}^{2})$ 
    & $\epsilon$ ($\mathbb{R}^{3}$) \\
    \hline
    0  & $7.9753\times 10^{-4}$ & $1.3022\times 10^{-3}$ & 
         $6.2537\times 10^{-4}$ \\
    1  & $2.0001\times 10^{-2}$ & $3.0382\times 10^{-2}$ & 
         $1.2038\times 10^{-2}$ \\
    2  & $5.2761\times 10^{-2}$ & $8.3471\times 10^{-2}$ & 
         $3.1985\times 10^{-2}$ \\
    4  & $1.1366\times 10^{-1}$ & $1.8888\times 10^{-1}$ & 
         $7.1597\times 10^{-2}$ \\
    10 & $2.3711\times 10^{-1}$ & $4.1171\times 10^{-1}$ & 
         $1.5639\times 10^{-1}$
  \end{tabular}
  \caption{Values of $\epsilon$ and $n$ used in
  Figures~\ref{fig:r1_oscillator_error2},
  \ref{fig:r2_oscillator_error2} and \ref{fig:r3_oscillator_error2}.
  With the presented initial conditions, each pairing gives a
  regularized kernel that results in the same modelling error,
  equation~\eqref{eqn:modelling_error}, of $4.89\times10^{-6}$.}
  \label{tbl:n_eps_pairs3} 
\end{SCtable}

\subsection{Harmonic oscillator in $\mathbb{R}^{1}$}
\label{sec:r1}

\begin{figure}[htbp]
  \centering
  \begin{minipage}{0.48\textwidth}
    \includegraphics[scale=1]{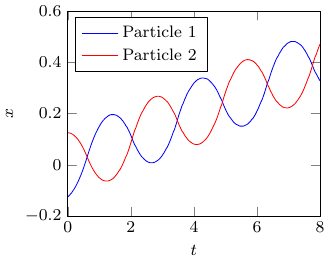}
  \end{minipage}
  \begin{minipage}{0.48\textwidth}
    \includegraphics[scale=1]{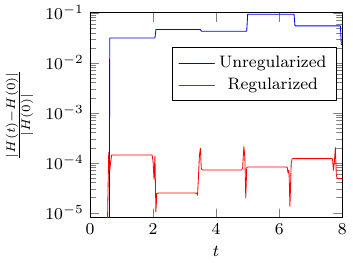}
  \end{minipage}
  \caption{\label{fig:r1_positions} The positions of the particles
    (left) and the errors in the Hamiltonian (right) when using the
    unregularized potential and regularized potential in
    $\mathbb{R}^{1}$.  Time stepping was done with a fourth-order
    symplectic integrator.  In this example, the time step size is
    $3.125\times 10^{-2}$ and the smoothing regularization parameter
    pairing is $(\epsilon,n) = (5.6755 \times 10^{-2},10)$.}
\end{figure}

The unregularized system is simply
\begin{align*}
  \ddot{x}_{j} = \frac{1}{2}\sum_{k \neq j} \left\{
    \begin{array}{rl}
      -1  & \quad x_{j} > x_{k}, \\
       1  & \quad x_{j} < x_{k}.
    \end{array}
  \right. 
\end{align*}
We consider two particles initially located at -0.125 and 0.125.  To
break the symmetry of the problem, which can cause errors to cancel, we
set the initial velocity of the left particle to be 0.1 and the right
particle to be 0.  Applying a fourth-order symplectic integrator, we see
in Figure~\ref{fig:r1_positions} that whenever the particles cross,
there is a jump in the error of the Hamiltonian.  This jump is caused
by the lack of regularity of the derivative of the Green's function.
If there are only a few particles, it is possible to exactly fix the
jump in the Hamiltonian by using an adaptive time step size.  However,
this strategy is not practical for many particles, or in higher
dimensions.  As an alternative, we replace the singular kernel with a
regularized kernel.  While the jumps are still present when using a
regularized kernel because of large derivatives, they are much smaller
than those for the unregularized system.

\begin{figure}[htbp]
    \centering
    \includegraphics[scale=1]{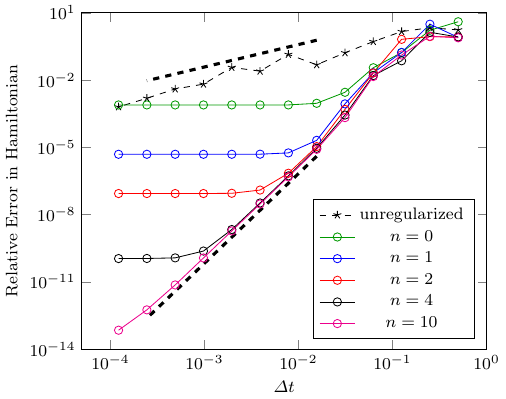}
  \caption{\label{fig:r1_oscillator_error} The error in the Hamiltonian
  arising from using kernel regularization in $\mathbb{R}^{1}$ with
  pairings of $(\epsilon,n)$ in Table~\ref{tbl:n_eps_pairs2}.  While
  fourth-order convergence is eventually achieved for all $(\epsilon,n)$
  pairs, the error eventually plateaus due to the modelling error of
  solving a regularized system.  Although the $(\epsilon,n)$ pairs have
  the same global regularization error, smaller modelling errors can be
  achieved when the higher-order kernels are used.  The dashed black
  lines correspond to first- and fourth-order convergence.}
\end{figure}

While the use of regularized kernels reduces the size of the jumps in
the Hamiltonian error, it does introduce a modelling error.  In
Figure~\ref{fig:r1_oscillator_error}, we plot the error in the
Hamiltonian for six different regularizations: the unregularized
kernel, and kernels regularized with the $(\epsilon,n)$ pairings in
Table~\ref{tbl:n_eps_pairs2}.  These pairings are specifically chosen
because each $(\epsilon,n)$ pair has the same global regularization
error.  For large $\Delta t$, there is no benefit in using the
high-order regularized kernels.  However, if smaller errors need to be
achieved, then it is favorable to use a regularized kernel.
Furthermore, we see that our new regularized kernels, $n>0$, achieve
smaller modelling errors, even though they have the same global
smoothing error defined in equation~\eqref{eqn:total_err}.

\begin{figure}[htbp]
  \begin{minipage}{0.7\textwidth}
    \centering
    \includegraphics[scale=1]{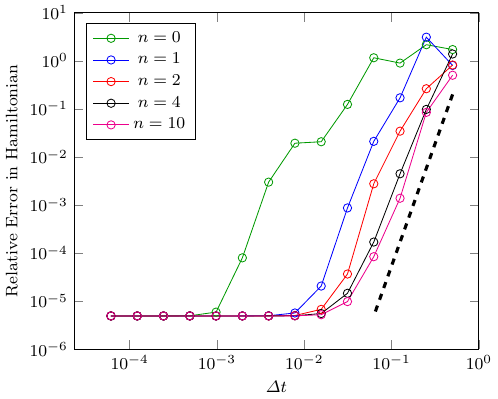}
  \end{minipage}\hfill
  \begin{minipage}{0.29\textwidth}
  \caption{\label{fig:r1_oscillator_error2} The error in the
    Hamiltonian arising from using kernel regularization in
    $\mathbb{R}^{1}$ with pairings of $(\epsilon,n)$ in
    Table~\ref{tbl:n_eps_pairs3}.  Higher-order
    kernels have smoother transition regions, allowing for high-order
    convergence with larger time step sizes.  The dashed black line
    corresponds to fourth-order convergence.}
  \end{minipage}
\end{figure}

It could be argued that the regularization error can be simply
decreased by taking a smaller value for $\epsilon$ while keeping $n$
fixed.  However, if $\epsilon$ is decreased, the derivative of the
regularized kernel increases at the origin, and the result is a smaller
asymptotic region for fourth-order convergence.  In
Figure~\ref{fig:r1_oscillator_error2}, we compare the error in the
Hamiltonian for $(\epsilon,n)$ pairings that all have a modelling error
of $4.89 \times 10^{-6}$.  Using larger values of $n$ 
results in larger regions of fourth-order convergence.  The trade-off
is the increased computational complexity for evaluating higher-order kernels.

\begin{figure}[htbp]
  \begin{minipage}{0.48\textwidth}
    \includegraphics[scale=1]{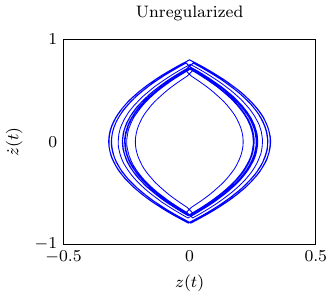}
  \end{minipage}
  \begin{minipage}{0.48\textwidth}
    \includegraphics[scale=1]{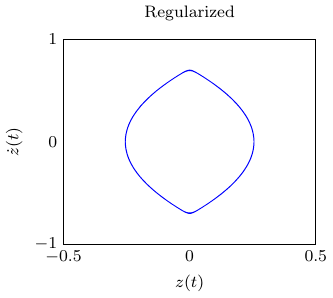}
  \end{minipage}
  \caption{\label{fig:phase_plane}The phase field of the $z(t) =
    x_{1}(t) - x_{2}(t)$.  The regularized potentials maintain the
    periodic solution of the problem for coarse time steps.}
\end{figure}

Finally, we examine the phase plane of the variable $z(t) = x_{1}(t) -
x_{2}(t)$.  We increase the time horizon from $T=8$ to $T=400$ and keep
the time step size fixed at $3.125 \times 10^{-2}$.  In
Figure~\ref{fig:phase_plane}, we plot the position $z(t)$ versus the
velocity $\dot{z}(t)$ resulting from the unregularized potential and the
regularized potential with $n=10, \epsilon=5.6755\times10^{-2}$; the
phase portrait with the other $(\epsilon,n)$ pairings in
Table~\ref{tbl:n_eps_pairs2} are indistinguishable in the eyeball norm
from the $n=10$ kernel.  The phase portrait of the unregularized system
shows the effect of the truncation error due to a discrete time
integrator being used in conjunction with the singular kernel.  The
qualitative periodic nature of the oscillations are perturbed.  On the
other hand, the regularized kernels significantly reduce the truncation
error, better preserving the periodic nature of the orbits.  If one were
to compare the final solutions using the regularized kernels using the
($\epsilon,n$) pairings in Table~\ref{tbl:n_eps_pairs2}, the errors
would decrease with increasing $n$ (not shown).

\subsection{Harmonic oscillator in $\mathbb{R}^{2}$}
\begin{figure}[htbp]
  \centering
  \begin{minipage}{0.48\textwidth}
  \includegraphics[scale=1]{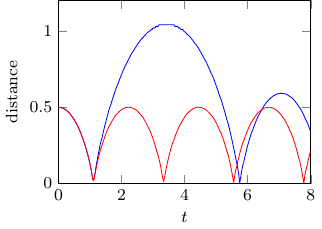}
  \end{minipage}
  \begin{minipage}{0.48\textwidth}  
  \includegraphics[scale=1]{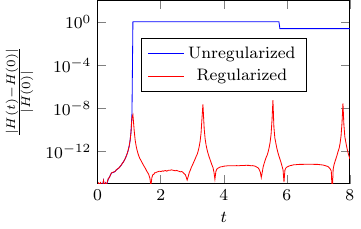}
  \end{minipage}
  \caption{\label{fig:r2_positions} The distance between the particles
  (left) and the errors in the Hamiltonian (right) when using the
  unregularized potential and regularized potential in
  $\mathbb{R}^{2}$.  Time stepping was done with a fourth-order
  symplectic integrator.  In this example, the time step size is $4.88
  \times 10^{-4}$ and the smoothing regularization parameter pairing is
  $(\epsilon,n) = (3.6132 \times 10^{-2},10)$.}
\end{figure}

Two particles are initially placed at $(-0.25,0)$ and $(0.25,0)$
with initial velocities $(0,10^{-3})$ and $(0,0)$ respectively.  The
unregularized system is
\begin{align*}
  \ddot{\xx}_{j} = \frac{1}{2\pi}\sum_{k \neq j} 
  \frac{\xx_{j} - \xx_{k}}{|\xx_{j} - \xx_{k}|^{2}}.
\end{align*}
These initial conditions were chosen so that the particles come close
to each other without actually passing through each other.  We expect
the close proximity of the particles to each other to delay the
fourth-order convergence if the unregularized system is solved.  In
Figure~\ref{fig:r2_positions}, the distance between the two particles
and the error in the Hamiltonian as a function of time is plotted.  
When the particles are close to each other, there is a jump
in the error of the Hamiltonian.

\begin{figure}[htbp]
  \centering
  \includegraphics[scale=1]{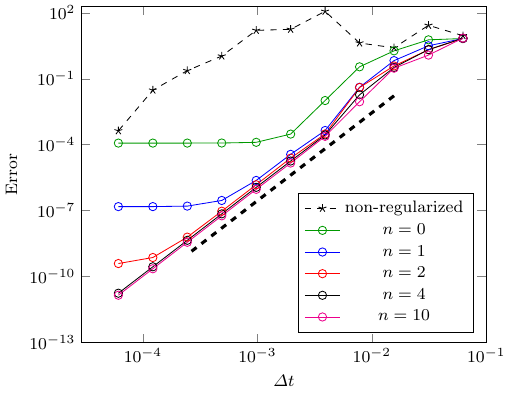}
  \caption{\label{fig:r2_oscillator_error} The error in the
    Hamiltonian arising from using $(\epsilon,n)$ pairings in
    Table~\ref{tbl:n_eps_pairs2} for the oscillator in
    $\mathbb{R}^{2}$.  Fourth-order convergence is achieved for both
    the unregularized kernel and the regularized kernels.  However,
    for the unregularized kernel, smaller time steps are required to
    enter this asymptotic regime. For the regularized systems, the
    error eventually plateaus when the modelling error dominates.  By
    using larger values of $n$, smaller modelling errors can be
    achieved.  The dashed black line corresponds to fourth-order
    convergence.}
\end{figure}

\begin{figure}[htbp]
  \begin{minipage}{0.7\textwidth}
    \centering
    \includegraphics[scale=1]{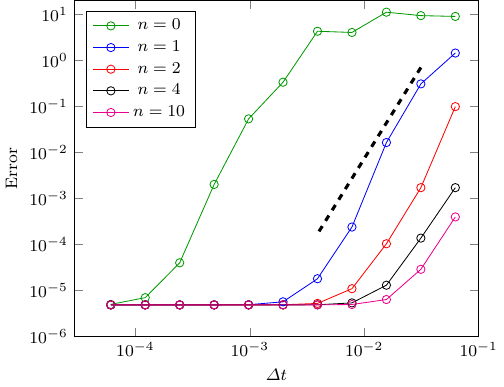}
  \end{minipage}\hfill
  \begin{minipage}{0.29\textwidth}
  \caption{\label{fig:r2_oscillator_error2} The error in the
  Hamiltonian arising from using kernel regularization in
  $\mathbb{R}^{2}$ with pairings of $(\epsilon,n)$ that have the same
  modelling error (as opposed to the smoothing error) of $4.89 \times
  10^{-6}$.  Larger values of $n$ achieve desired accuracies with
  larger time step sizes.  The dashed black line corresponds to
  fourth-order convergence.
  }
  \end{minipage}
\end{figure}

Repeating the numerical experiments from Section~\ref{sec:r1},
Figure~\ref{fig:r2_oscillator_error} shows the convergence behavior of
the regularized dynamical system if the $(\epsilon,n)$ pairings from
Table~\ref{tbl:n_eps_pairs2} (fixed global smoothing error) are used.
Similar observations to the experiment in $\mathbb{R}^1$ can be
observed: the symplectic integrator achieves fourth-order accuracy for
each regularized system for large $\Delta t$, until the modeling error
dominates.  If the $(\epsilon,n)$ pairings from
Table~\ref{tbl:n_eps_pairs3} (fixed modelling error) are used, the
higher-order regularized kernels exhibit fourth-order convergence for
larger $\Delta t$ (c.f.~Figure~\ref{fig:r2_oscillator_error2}).
%
%

\subsection{Harmonic oscillator in $\mathbb{R}^{3}$}

\begin{figure}[htbp]
  \begin{minipage}{0.7\textwidth}
    \centering
    \includegraphics[scale=1]{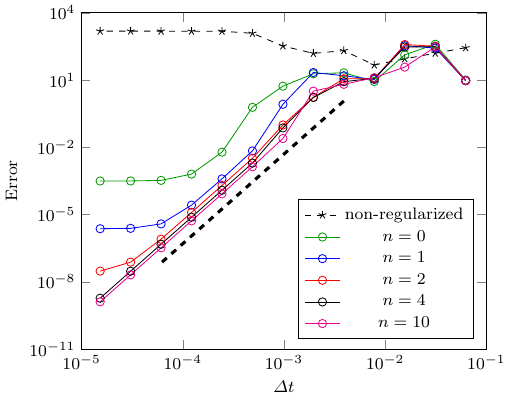}
  \end{minipage}\hfill
  \begin{minipage}{0.29\textwidth}
  \caption{\label{fig:r3_oscillator_error} The error in the Hamiltonian
  arising from using kernel regularization in $\mathbb{R}^{3}$.  While
  fourth-order convergence is achieved for all values of $n$, the error
  eventually plateaus.  By using larger values of $n$, smaller errors
  can be achieved.  The dashed black line corresponds to fourth-order
  convergence.}
  \end{minipage}
\end{figure}

\begin{figure}[htbp]
  \begin{minipage}{0.7\textwidth}
    \centering
    \includegraphics[scale=1]{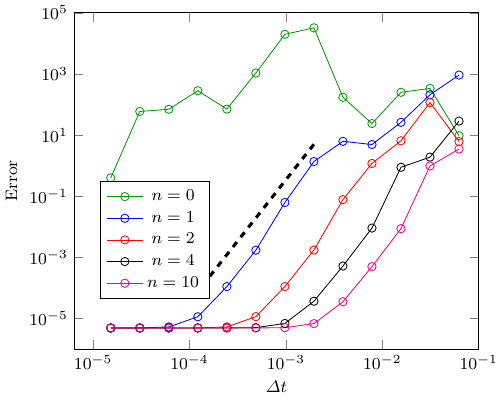}
  \end{minipage}\hfill
  \begin{minipage}{0.29\textwidth}
  \caption{\label{fig:r3_oscillator_error2} The error in the Hamiltonian
  arising from using kernel regularization in $\mathbb{R}^{3}$ with
  pairings of $(\epsilon,n)$ that have the same regularization error of
  $4.89 \times 10^{-6}$.  Larger values of $n$ achieve desired
  accuracies with larger time step sizes.  The dashed black line
  corresponds to fourth-order convergence.
  }
  \end{minipage}
\end{figure}

Similarly, we place 2 particles at $(-0.1,0,0)$ and $(0.1,0,0)$ with
initial velocities of $(0,10^{-3},0)$ and $\vec{0}$ respectively.  The
unregularized system is
\begin{align*}
  \ddot{\xx}_{j} = \frac{1}{4\pi}\sum_{k \neq j} 
  \frac{\xx_{j} - \xx_{k}}{|\xx_{j} - \xx_{k}|^{3}} .
\end{align*}
As before, we expect that the singularity will reduce the order of
accuracy of the fourth-order symplectic integrator.  In fact, since the
singularity is even stronger than in $\mathbb{R}^{2}$, we see that the
unregularized system does not even obtain convergence for the reported
values of  $\Delta t$ (Figure~\ref{fig:r3_oscillator_error}).  By
introducing a regularization, fourth-order convergence is observed.
Again, as before, we observe that if the global smoothing error is kept
constant, than larger values of $n$ reduce the modelling error, and we
are able to achieve more accurate results.  As in $\mathbb{R}^{1}$, the
benefit of using large values of $n$ is illustrated in
Figure~\ref{fig:r3_oscillator_error2}.  The modelling error is fixed at
$4.89 \times 10^{-6}$, and the error in the Hamiltonian is plotted for
the different $(\epsilon,n)$ pairings.

\begin{figure}[htbp]
  \begin{minipage}{0.48\textwidth}
    \includegraphics[scale=1]{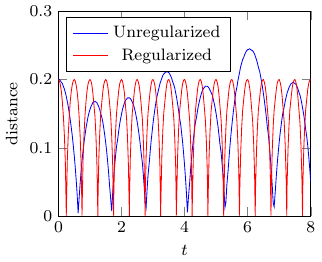}
  \end{minipage}
  \begin{minipage}{0.48\textwidth}      
    \includegraphics[scale=1]{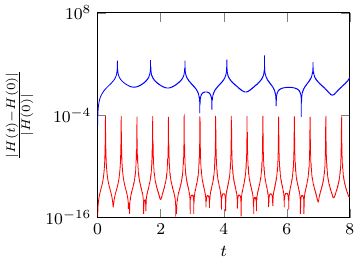}
  \end{minipage}
  \caption{\label{fig:r3_positions} The distance between the particles
  (left) and the errors in the Hamiltonian (right) when using the
  unregularized potential and regularized potential in
  $\mathbb{R}^{3}$.  Time stepping was done with a fourth-order
  symplectic integrator.  In this example, the time step size is
  $2.44\times 10^{-4}$ and the smoothing regularization parameter
  pairing is $(\epsilon,n) = (2.8378 \times 10^{-2},10)$.}
\end{figure}

\subsection{Periodic orbit in $\mathbb{R}^3$}
\begin{figure}[htbp]
  \centering
  \begin{minipage}{0.45\textwidth}
    \begin{center}
      \includegraphics[scale=1]{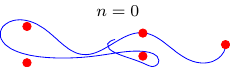} \\
      \vspace{10pt}
      \includegraphics[scale=1]{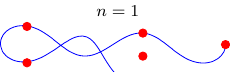}
    \end{center}
  \end{minipage}
  \begin{minipage}{0.45\textwidth}
    \begin{center}
      \includegraphics[scale=1]{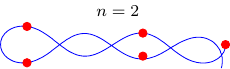} \\
      \vspace{10pt}
      \includegraphics[scale=1]{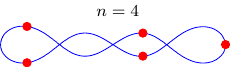}
    \end{center}
  \end{minipage}
  \caption{\label{fig:periodic_orbits}The initial location of the five
  particles and the trajectory that one of the particles follows over
  one period.  If $n$ is too small, the correct orbit can not be
  achieved.  However, with $n=4$, the error in periodicity is
  $10^{-2}$.  The orbit with $n=10$ looks similar to the orbit with
  $n=4$.}
\end{figure}

We consider five particles in the $z=0$ plane of $\mathbb{R}^{3}$.  We
change the sign of $G(r)$ so that the Hamiltonian system corresponds to
motion due to the gravitational potential.  By setting the mass of each
particle to $\sqrt{0.2}$, so that $w_{jk} = 2$, and using the initial
condition 
\begin{align*}
  \xx(0) = \left(
  \begin{array}{c@{\hspace{20pt}}c}
   +3.315332 \times 10^{-1} & 
   0 \\
   +8.795500 \times 10^{-2} &
   -3.394340 \times 10^{-2} \\
   -2.537216 \times 10^{-1} &
   -5.353020 \times 10^{-2} \\
   -2.537216 \times 10^{-1} &
   +5.353020 \times 10^{-2} \\
   +8.795500 \times 10^{-2} &
   +3.394340 \times 10^{-2}
  \end{array}
  \right), \\
  \dot{\xx}(0) = \left(
  \begin{array}{c@{\hspace{20pt}}c}
   0 &
   -5.937860 \times 10^{-1} \\
   +1.822785 \times 10^{0} &
   +1.282480 \times 10^{-1} \\
   +1.271564 \times 10^{0} &
   +1.686450 \times 10^{-1} \\
   -1.271564 \times 10^{0} &
   +1.686450 \times 10^{-1} \\
   -1.822785 \times 10^{0} &
   +1.282480 \times 10^{-1}
  \end{array}
  \right),
\end{align*}
the dynamics should result in a periodic orbit with period
$T=2\pi/5$~\cite{simo2001}.  Using a fourth-order symplectic integrator
with $10^6$ time steps, the non-regularized system does not give a
periodic orbit,  due to the singularity of the kernel.  Using the
$(\epsilon,n)$ pairings from Table~\ref{tbl:n_eps_pairs2} (i.e., fixed
global smoothing error), the orbits of the regularized system are shown
in Figure~\ref{fig:periodic_orbits}.  For the $n=0, 1, 2$ kernels, the
modelling error dominates in the regularized system, resulting in
orbits that are qualitatively different from the expected periodic
orbit.  The $n=4$ and $n=10$ kernels result in regularized systems that
give qualitatively correct periodic orbits.  

The largest possible time step that keeps the relative error in
periodicity bounded at $10^{-2}$, i.e. $\|\xx(T)-\xx(0)\| <
10^{-2}\|\xx(0)\|$, is computed for the $n=4$ and $n=10$ systems.  The
results are reported in Table~\ref{tbl:fiveparticles}.  A larger time
step can be used for the $n=10$ kernels.
\begin{table}[htbp]
\centering
\begin{tabular}{ccccc}
$n$ & $\Delta t$ & Period Error & Hamiltonian Error & Modelling Error \\
\hline
4  & $1.58 \times 10^{-3}$ & $9.99 \times 10^{-3}$ &
     $3.65 \times 10^{-7}$ & $2.50 \times 10^{-7}$ \\
10 & $2.29 \times 10^{-3}$ & $9.98 \times 10^{-3}$ & 
     $5.13 \times 10^{-7}$ & $9.87 \times 10^{-11}$
\end{tabular}
\caption{\label{tbl:fiveparticles}A summary of the different errors
  when considering five particles that form a periodic orbit.  We see
  that using $n=10$ allows for a larger time step size than $n=4$.}
\end{table}

\section{Conclusions}
\label{sec:conclusions}
In this paper, we derived a new family of regularized kernels, suitable
for simulating a Hamiltonian system that contains the fundamental
solution of Laplace's equation using high-order time integrators.
These high-order kernels were obtained by a Taylor expansion of the
non-regularized kernel about $(r^2+\epsilon^2)$ in $\mathbb{R}^1,
\mathbb{R}^2$, and $\mathbb{R}^3$.  The analysis shows that the
regularized kernels, $G^{\epsilon,n}(r)$, converge to the fundamental
solution of Laplace's equation as $n\to\infty$ for any $\epsilon>0$,
and as $\epsilon\to0$ for any $n\ge 0$.  In addition to the derivation
and validation of the high-order kernels, error bounds for the
regularized solution were derived.

We have shown that these regularizations can reduce the error in the
far field without introducing sharp derivatives near the singularity.
This is particularly useful when applying high-order time stepping
methods to a Hamiltonian system such as a harmonic oscillator.  In
particular, high-order regularized kernels (with identical global
smoothing error) can reduce the modelling error of the regularized
system. Alternatively, if one chooses regularizations that
give rise to similar modelling errors, high-order accuracy can be
achieved for larger time step sizes using these high-order kernels.

Future work includes using these high-order kernels within treecode
algorithms for approximating the electric field arising from a charged
particle
system~\cite{christlieb:tas04,christlieb:ckv03IEEE,salmon-warren94}---this
is necessary when simulating a large number of particles.  In a
treecode computation, a Taylor series expansion of the regularized
potential is needed for the computation of the cluster-particle
interaction.  While explicit formulas for the Taylor coefficients of
the high-order regularized kernels can be derived and evaluated, the
computational complexity might be prohibitive.  The authors anticipate
that a recurrence relation to evaluate the Taylor coefficients of the
high-order regularized kernels can be recovered.  Certainly, the
recurrence relation for the $n=0$ kernels are
available~\cite{lindsay-krasny01}.  Alternatively, a kernel-independent
fast multipole method~\cite{ros-ols2015} can be used to significantly
reduce the number of computations.  With these fast algorithms, more
complicated simulations such as vortex
motions~\cite{krasny1986dpv,krasny1987computation} can be
investigated.  High-order regularized kernels can likely also be
formulated in a similar fashion for for the screened Coulomb
potential~\cite{Li20093858} or Winckelmans--Leonard
kernel~\cite{MR2520287}.


\section*{Acknowledgments}

The authors would like to thank Robert Krasny, Keith Cartwright, John
Verboncoeur, John Luginsland, Matthew Bettencourt, and Andrew
Greenwood for their insightful discussions regarding this work, as
well as anonymous referees who have made valuable suggestions to
improve the presentation of this manuscript.


\end{document}